\documentclass[12pt]{article}
\usepackage{amsmath}

\newtheorem{theorem}{Theorem}

\newtheorem{lemma}{Lemma}

\newtheorem{proposition}{Proposition}
\newtheorem{remark}[theorem]{Remark}

\numberwithin{equation}{section}

\begin{document}

\title{On spectral representations  of \\ tensor random fields on the sphere\thanks{
Partly supported by the Welsh Institute of Mathematics and
Computational Sciences and the Commission of the European
Communities grant PIRSES-GA-2008-230804 within the programme
`Marie Curie Actions'. }}
\author{N. Leonenko \thanks{%
Cardiff School of Mathematics, Cardiff University, Senghennydd
Road, Cardiff,CF24 4AG, UK, Email:LeonenkoN@Cardiff.ac.uk} \ \ \
and \ \ \ L.
Sakhno \thanks{%
Dept. of Probability Theory and Mathematical Statistics, Kyiv
National Taras Shevchenko University, Ukraine,
Email:lms@univ.kiev.ua}}


\maketitle

\begin{abstract}
\noindent We study the representations of tensor random fields on
the sphere basing on the theory of representations of the rotation
group. Introducing specific components of a tensor field and
imposing the conditions of weak isotropy and mean square
continuity, we derive their spectral decompositions in terms of
generalized spherical functions. The properties of random
coefficients of the decompositions are characterized, including
such an important question as conditions of Gaussianity.



\end{abstract}



\section{Introduction}

This paper was inspired and motivated by recent advances in the study of
cosmic microwave background (CMB) radiation, which is currently a crucial
topic of investigation in cosmology. Amazingly, the CMB anisotropies contain
a wealth of important cosmological information and provide the main testing
ground for the theories of early universe, with the possibility of precise
determination of fundamental cosmological parameters. Note that this area
has raised a number of issues, theoretical and practical, both in physics
and mathematics. Nowadays, with the flood of high quality data available
from several satellite missions and expected from the future ones, the focus
of CMB research is on developing theoretical tools for data interpretation
and analysis.

The CMB radiation is characterized by four Stokes parameters
\cite{Chandrasekhar}: the intensity $%
I$ (or, equivalently, the temperature $T$) and the other three parameters $Q$%
, $U$ and $V$, which define the polarization state, namely, linear and
circular polarization (however, the fourth parameter $V$ describing the
circular polarization is not necessary in the standard cosmological models).

Much more attention has been paid to the study of temperature fluctuations
in cosmic microwave background, great deal of experimental activity has been
accompanied by thorough numerical and analytic work. In particular, rigorous
mathematical basis has been provided through the excellent series of papers
\cite{BaldiMarinucci}, \cite{BaldiMarinucci2}, \cite{Marinucci2005}--\cite
{MarinucciPiccioni}.

From the mathematical point of view, the temperature fluctuations are
considered as (a realization of) a random field on the unit sphere. One
natural tool for statistical analysis of random fields is a spectral
decomposition, which for scalar random fields on the sphere is a
decomposition in the series of spherical harmonics (see, e.g., \cite
{Leonenko}, \cite{Yadrenko}). These decompositions have been known and used
in different applied areas for quite a long time, however recent
contributions in \cite{BaldiMarinucci}, \cite{BaldiMarinucci2}, \cite
{Marinucci2005}--\cite{MarinucciPiccioni} shed more light on many important
practical issues such as characterization of spherical harmonics
coefficients and their asymptotic behavior, characterization of spectrum and
bispectrum, testing for Gaussianity on spherical random fields and some
other. Higher order spectra for spherical random fields have been analyzed
in \cite{MarinucciPeccati}.

It has been long recognized that polarization in CMB provides a
source of information on cosmological parameters complementary to
temperature fluctuations. This significant additional information
is necessary to create a complete picture of anisotropies used to
constrain cosmological models, to improve the accuracy in
determining cosmological parameters, and also to study relic
gravitational waves, which is one of the biggest challenges today
(see, e.g., \cite{Zhao}, \cite{Baskaran}). However, polarization
anisotropies in CMB are detected at much smaller scales than
temperature fluctuations and this limitation has been overcome
only recently with improved sensitivity and resolution of CMB
experiments. First polarization measurements were obtained by WMAP
experiment (see, e.g., \cite{Page}, http://map.gsfc.nasa.gov/) and
further data are expected from Plank surveyor satellite (see
http://www.sciops.esa.int/) prompting many groups of scientists to
focus on polarization research. Last few years of intensive
studies in this rapidly expanding area resulted in developing
physical and mathematical formalism to describe, characterize and
analyze the polarization counterpart in CMB.

From the mathematical point of view, considering polarization we deal again
with random fields on the sphere, however, polarization Stokes parameters $Q$
and $U$ are not scalars but rather the components of a symmetric trace-free
rank 2 tensor, or, equivalently, the combinations $Q\pm iU$ are spin $\pm 2$
quantities. Therefore, for their mathematical description one needs to
consider tensor-valued random fields (or random fields with spin-weighted
quantities) on the sphere and their spectral decomposition and statistical
analysis.

In order to study polarization anisotropies, analytic methods have
been formulated in physical literature based on expansions of
tensor-valued functions on the sphere in the series of appropriate
tensor spherical harmonics, or expansions in spin-weighted
spherical harmonics (this was done initially by two group of
scientists in \cite{Kamionkowski} and \cite {Zaldarriaga}
respectively, and numerous further developments have been
elaborated by now). The main reference for the mathematical
formalism used in these studies comes from non-probabilistic
literature, in particular, the tools for the mentioned above
expansions are provided by the group representation theory,
harmonic analysis on the group of rotations and the Peter-Weyl
theorem. However, more caution and accuracy are needed when
considering expansions for random functions, in particular, one
needs to precise in what sense and under what conditions the
corresponding series with random coefficients converge. These
important issues are missing in physical literature, or at least
they are not rigorously and clearly pronounced. It is also
important to reveal in detail the nature and features of the
random coefficients of the series expansions. Note that rigorous
probabilistic framework is necessary to proceed with statistical
inference methods.

Only recently the rigorous mathematical and probabilistic basis to
underpin used in physical literature techniques appeared in
several papers: \cite{GellerMarinucci2008},
\cite{GellerLanMarinucci2009}, \cite {Malyarenko2009}. In these
papers the theory for random sections of fiber bundles over the
sphere is developed with cosmological applications in mind.

In this paper we present an approach to derive spectral
decompositions of tensor-valued random fields on the sphere from
more group-theoretical point of view. We introduce specific
components of a tensor field and under the conditions of weak
isotropy and mean square continuity define the representation of
the rotation group by rotational shift operators acting on the
Hilbert spaces generated by these components. Then spectral
decompositions of the components are obtained as the sum of
projections on the irreducible representation spaces, which
expressed in terms of generalized spherical functions. We
characterize the properties of random coefficients of the
expansions obtained, including such an important question as
conditions of Gaussianity.

Note that problems involving spectral analysis for random fields on the
sphere arise in various applied areas such as geodesy, geophysics, planetary
sciences, astronomy, cosmology, medical imaging etc.

The organization of the paper is as follows. We start in Section 2 with the
brief overview of the representation theory of the rotation group, to that
extend which will be needed in the further exposition. In Section 3 we
review different approaches used to derive the spectral representation for
scalar random fields on the sphere, with the particular attention paid to
the approach based on representation of the rotation group with
representation space, which is the Hilbert space generated by the random
field. This last approach will be applied in Section 4 to derive spectral
representations of vector and tensor fields.

\section{Elements of the theory of the rotation group $SO(3)$ and its
representations}

We summarize briefly some important facts on the group $SO(3)$, which will
be used in the following sections. In our exposition we prefer to follow to
the classical books \cite{Gel}, \cite{VK}, we will also indicate the
connection to some other definitions and notions used in the literature, in
particular, we will refer to yet another classical textbook \cite{Varsh}.
However, many other excellent sources on the topic are available.

The group $SO(3)$ of special orthogonal transformations in $R^{3}$ consists
of all rotations in $R^{3}$ about a fixed origin. This group has a
convenient realization as the group of all $3\times 3$ real matrices $A$
such that $A^{t}A=I$ and $detA=1$.

A number of parametric representations of rotations can be introduced. One
important parametric form for the group of rotation uses the so-called Euler
angles. Note that there exist different conventions how to define these
angles. One of such conventions (e.g., \cite{Gel}) is to define a rotation $%
g=g(\varphi_1, \theta, \varphi_2)$ as the product of three successive
rotations: a rotation $g_{\varphi_1}$ by an angle $\varphi_1$ around the $z$
axis, then a rotation $g_{\theta}$ by an angle $\theta$ around the new $x$
axis, then a rotation $g_{\varphi_2}$ by an angle $\varphi_2$ around the new
$z$ axis. 
Correspondingly, the elements of the matrix of the rotation can be expressed
explicitly in terms of the Euler angles $(\varphi_1, \theta, \varphi_2)$
(see, e.g., \cite{Gel}).

Let $V$ be a Hilbert space. A representation of a group $G$ with a
representation space $V$ is a homomorphism $T: g \to T(g)$ of $G$ into the
space of bounded linear operators on $V$. Thus, the mapping $T$ satisfies
the conditions:
$$
T(g_1)T(g_2)=T(g_1g_2), \ \ \ T(e)=E,
$$
where $e$ is an identity element of $G$, $E$ is an identity operator. Note,
that with each operator $T(g)$ one can associate its matrix $\{t_{ij}(g)\}$
(defined with respect to an orthonormal basis in V).

A representation $T$ is said to be irreducible if there does not exist a
proper subspace $W$ of $V$, which is invariant under $T$. Otherwise, $T$ is
reducible. The representation $T$ is said to be unitary if the operators $%
T(g)$, $g\in G$, are unitary with respect to the scalar product defined on $%
V $, i.e., for all $x, y \in V$, $g\in G$ we have
$$
(T(g)x, T(g)y) =(x,y).
$$

We formulate now the basic results on the resolution of a representation
into irreducible representations (see, \cite{Gel}).

\begin{proposition}
Let $T: g\to T(g)$ be an unitary representation of the group of rotations in
a (separable) Hilbert space $V$. Then there exist mutually orthogonal
finite-dimensional subspaces $V_1, V_2, \ldots$, invariant with respect to $%
T $, in each of which the representation $T$ is irreducible, and the space $%
V $ is the orthogonal sum of these subspaces $V_i$. This means that every $%
x\in V $ can be expresses as a convergent series $x=\sum x_i$, $x_i\in V_i$,
where convergence is meant with respect to the norm generated by the scalar
product defined on $V$.
\end{proposition}

\begin{proposition}
Each irreducible representation of the group of rotations is defined by a
number $l$, called the weight of representation, the corresponding invariant
representation space $D_{l}$ has dimension $2l+1$. The matrix of the
representation (of weight $l$) corresponding to an arbitrary rotation $%
g=g(\varphi _{1},\theta ,\varphi _{2})$ has, in the canonical basis, the
form
$$
T^{l}(g)=\{T_{mn}^{l}(\varphi _{1},\theta ,\varphi _{2})\}_{m,n=-l,...,l},
$$
where the elements are given by the following functions
\begin{equation}
T_{mn}^{l}(\varphi _{1},\theta ,\varphi _{2})=e^{-im\varphi
_{1}}P_{mn}^{l}(\cos \theta )e^{-in\varphi _{2}}.  \label{tlmn}
\end{equation}
\end{proposition}

Note that the matrix representation (\ref{tlmn}) is referred to the fixed
canonical basis in the space $D_l$, which is formed by $2l+1$ orthogonal
vectors (of dimension $2l+1$), each being the eigenvector with the
eigenvalue $e^{-in\varphi}$, $n=-l,\ldots,l$, for the rotation around $z$
axis by an angle $\varphi$.

The matrix elements $T^l_{mn}(\varphi_1, \theta, \varphi_2)$ are called the
generalized spherical functions of the order $l$ (\cite{Gel}).

The functions $P^l_{mn}(z)$ appearing in the formula (\ref{tlmn}) can be
represented in different ways. We give here their differential
representation, which is also called the Rodrigues formula\footnote{%
The formula (\ref{plmn}) is actually called the Rodrigues formula in \cite
{VK} when $P^l_{mn}(z)$ is defined via (\ref{plmn}) with $i^{m-n}$ dropped,
instead $i^{m-n}$ is included in (\ref{tlmn}).}:
\begin{eqnarray}  \label{plmn}
P^l_{mn}(z) &=& \frac{(-1)^{l-m}i^{m-n}}{2^l}\left[ \frac{(l+n)!}{%
(l-m)!(l+m)!(l-n)!}\right]^{1/2}   \\
&\times& (1-z)^{-(n-m)/2}(1+z)^{-(n+m)/2}\frac{d^{l-n}}{dz^{l-n}}\left[%
(1-z)^{l-m}(1-z)^{l+m}\right]. \nonumber
\end{eqnarray}

Expressions for $P^l_{mn}(z)$ can be given in terms of hypergeometric
functions, trigonometric functions, there exists also their integral
representation. All this can be found, for example, in \cite{VK}, \cite
{Varsh}, where the connection of $P^l_{mn}(z)$ with classical orthogonal
polynomials and some other their properties are also reported.

In many cases one considers representations of a group $G$ by shift
operators in the linear space of functions defined on some space $X$, in
particular, $X$ can be equal to $G$ itself. These shift operators acting on
functions $f(g)$, $g\in G$, can be defined as
$$
T(g_{0})f(g)=f(g_{0}^{-1}g),
$$
which corresponds to the case when $G$ is considered as the group of left
shifts, or one can define the action of the operator $T(g_{0})$ as \
$$
T(g_{0})f(g)=f(gg_{0}),
$$
in the case when $G$ is considered as the group of right shifts.

Let $L_{2}(G)$ be the Hilbert space of all functions on the compact group $G$%
, which are square integrable with respect to the Haar measure $dg$ on $G$
(this measure is invariant for compact groups). Considering the case $%
G=SO(3) $, we have the space of functions
$$
f(g)=f(\varphi _{1},\theta ,\varphi _{2}),
$$
for which the following integral exists
\begin{equation}
\int_{G}|f(g)|^{2}dg=\frac{1}{8\pi ^{2}}\int \int \int \left| f(\varphi
_{1},\theta ,\varphi _{2})\right| ^{2}\sin \theta d\varphi _{1}d\theta
d\varphi _{2}<\infty ;
\end{equation}
the scalar product of $f_{1}(g)$ and $f_{2}(g)$ is defied as
\begin{equation}
(f_{1},f_{2})=\int f_{1}(g)\overline{f_{2}(g)}dg.
\end{equation}

The transformation
$$
T(g_{0})f(g)=f(gg_{0})
$$
forms a unitary representation in the space $L_{2}(G)$ (called the (right)
regular representation of the rotation group \cite{Gel}). The irreducible
representations into which it can be resolved are the representations in the
subspaces of generalized spherical functions $T_{mn}^{l}(\varphi _{1},\theta
,\varphi _{2})$ for a fixed $l$ and $m$. Thus, the resolution of this
representation into irreducible representations means that every function $%
f\in L_{2}(SO(3)) $ can be expanded as a series in the functions $%
T_{mn}^{l}(\varphi _{1},\theta ,\varphi _{2})$.

\begin{proposition}
The set of the generalized spherical functions $T^l_{mn}(\varphi_1, \theta,
\varphi_2)$ ($l$ being an integer) forms a complete orthogonal system in the
space of functions $L_2(SO(3))$.
\end{proposition}

Note that the above statement is a particular case of the Peter-Weyl
theorem, which is one of the most important results of the harmonic analysis
on compact groups (see, e.g., \cite{VK}). We will return to this theorem
later in this section and consider its stochastic version.

The multiplicative law of the group representation implies the following
rule according to which the generalized spherical functions are added.

\begin{proposition}
(Addition formula for the generalized spherical functions.)
\begin{equation}
T_{mn}^{l}(g_{1}g_{2})=\sum_{s=-l}^{l}T_{ms}^{l}(g_{1})T_{sn}^{l}(g_{2}).
\end{equation}
\end{proposition}

Using the unitarity of the representation matrices, the above formula can
also be written in the following form
\begin{equation}  \label{addition}
T_{mn}^{l}(g_{1}g_{2}^{-1})=\sum_{s=-l}^{l}T_{ms}^{l}(g_{1})\overline{%
T_{ns}^{l}(g_{2})}.
\end{equation}

\bigskip

Let us return to the definition of the Euler angles $(\varphi _{1},\theta
,\varphi _{2})$. Another way to introduce them, most commonly used in
physical literature, e.g., in quantum mechanics, is based on so-called $zyz$%
-convention about rotations (the introduced above is $zxz$-convention), when
a rotation is defined as to be produced in the following three steps: a
rotation $g_{\varphi _{1}}$ by an angle $\varphi _{1}$ around the $z$ axis,
then a rotation $g_{\theta }$ by an angle $\theta $ around the new $y$ axis,
then a rotation $g_{\varphi _{2}}$ by an angle $\varphi _{2}$ around the new
$z$ axis.

Correspondingly, within this approach, the matrix elements of the
irreducible unitary representation of the weight $l$ of the rotation group $%
SO(3)$ are represented by the Wigner $D$-functions $D^l_{mn}(\varphi_1,
\theta, \varphi_2)$ (see, \cite{Wigner}, \cite{BrinkSatchler}, \cite{Varsh}%
). Proposition 3 holds for the functions $D^l_{mn}$, that is, every function
from $L_2(SO(3))$ can be expanded in a series in terms of the Wigner $D$%
-functions.

Note that in quantum mechanics these functions appear to be the elements of
the matrix representation of the rotation operator in the basis formed by
the angular momentum eigenvectors. Wigner $D$-functions play an important
role in various fields of modern physics including nuclear and molecular
physics.

Both conventions about the definition of the Euler angles are closely
related and the same rotation $R$ can be achieved with a simple adjustment
of angles: $R_{zyz}(\alpha,\beta,\gamma)=R_{zxz}(\alpha+\pi/2,\beta,\gamma-%
\pi/2)$ (see, \cite{Varsh}).

Correspondingly, there exists a simple relation between the functions $%
T^l_{mn}$ and $D^l_{mn}$ (see, \cite{Varsh}):
\begin{equation}
D^l_{mn}(\varphi_1, \theta, \varphi_2)=(-i)^{n-m}T^l_{mn}(\varphi_1, \theta,
\varphi_2).
\end{equation}

However, an advantage of the Wigner $D$-functions is that their factorized
form is given in terms of the real Wigner $d$-functions $d_{mn}^{l}(\theta )$
and the complex exponentials:
\begin{equation}
D_{mn}^{l}(\varphi _{1},\theta ,\varphi _{2})=e^{-im\varphi
_{1}}d_{mn}^{l}(\cos \theta )e^{-in\varphi _{2}},  \label{dlmn}
\end{equation}
where $d_{mn}^{l}$ are equal to $P_{mn}^{l}$ with the factor $(-i)^{n-m}$
dropped (cf. (\ref{tlmn})-(\ref{plmn})).

\bigskip

As we have seen above the group representation theory (outlined here for a
particular case of the group $G=SO(3)$) provides a tool -- the Peter-Weyl
theorem -- for decomposition of $L_{2}(G)$ space (endowed with the Haar
measure) into an orthogonal sum of finite dimensional spaces and,
correspondingly, the decomposition of functions from $L_{2}(G)$ into the sum
of their projections on these spaces which can be represented in terms of
appropriate basis functions.

Our main interest is in extension of this result which allows for
decomposition of random functions. One can observe immediately that such a
decomposition can be derived for random functions $X(g)$, $g\in G$, whose
trajectories $g\rightarrow X(g)$ are $P$-a.s. square integrable w.r.t. the
Haar measure (that is, belong to $L_{2}(G)$), corresponding decomposition
must be understood in the $L_{2}(G)$-sense, $P$-a.s. However, more profound
results can be obtained with the assumption of isotropy, that is, invariance
in law of a random function $X(y)$, $y\in Y$, under the action of a group $G$
(with $G$ being a topological compact group acting on $Y$, $Y$ can be equal
to $G$ itself). In such a case one can formulate a stochastic version of the
Peter-Weyl theorem (see, \cite{PeccatiPycke}) and the decomposition of $%
L_{2}(G)$ translates easily into construction of spectral representations
for isotropic random functions. To be more precise, for the case of the
rotation group $SO(3)$ the following result can be stated as a consequence
of the stochastic Peter-Weyl theorem (see, \cite{MarinucciPeccati}, \cite
{PeccatiPycke}).

\begin{proposition}\label{prop5}
Let $X(g)$, $g\in G=SO(3)$, be a square integrable
$(L^{2}(P(d\omega ))$, strictly isotropic random field, that is,
finite-dimensional distributions of $$\{X(g_{1}), \ldots ,
X(g_{k})\} \mbox{\ \ and\ \ } \{X(gg_{1}), \ldots, X(gg_{k})\}$$
are the same $\forall k$, $\forall g,$ $g_{1},\ldots ,g_{k}\in G$.
Then
\begin{equation}  \label{spectralg}
X(g)=X(\varphi ,\theta ,\psi )=\sum_{l}\sum_{m,n}a_{lmn}\sqrt{\frac{2l+1}{%
8\pi ^{2}}}D_{mn}^{l}(\varphi ,\theta ,\psi ),
\end{equation}
both in $L^{2}(P(d\omega )\times dg)$ and pointwise in $L^{2}(P(d\omega ))$,
where $dg$ is the Haar measure on $SO(3),$ and
\begin{equation}  \label{almn}
a_{lmn}=\int_{G}X(g)\sqrt{\frac{2l+1}{8\pi ^{2}}}\overline{D_{mn}^{l}}(g)dg.
\end{equation}
\end{proposition}

\bigskip

\begin{remark} {\rm The functions $D_{mn}^{l}(g)$ (as well as $T_{mn}^{l}(g)$)
are orthogonal but not orthonormal. The orthogonality relation for $%
D_{mn}^{l}(g)$ reads
$$
\int_{G}{D_{mn}^{l}}(g)\overline{D_{m^{\prime}n^{\prime}}^{l^{\prime}}}(g)dg=%
\frac{8\pi ^{2}}{2l+1}\delta_{ll^{\prime}}\delta_{mm^{\prime}}\delta_{nn^{%
\prime}}.
$$
This explains the appearance of the factor $\sqrt{\frac{2l+1}{8\pi
^{2}}}$ in (\ref{spectralg})-(\ref{almn}).} \end{remark}


\begin{remark} {\rm In what follows we will refer to the above result,
however we will also consider the representation of $SO(3)$ via
(rotational) shift operators defined on the Hilbert space
generated by a random field,
which gives rise to alternative derivation of the spectral representation (%
\ref{spectralg}) and allows to treat the coefficients from the
different point of view.}
\end{remark}

\section{Spectral representation of scalar fields}

In this section we review different approaches used to derive the spectral
representation for scalar random fields on the sphere.

In what follows we will suppose that we are given the probability space $%
(\Omega, \mathcal{F}, P)$.

Let $L^2(P(d\omega))$ be the Hilbert space of random variables $Y$ such that
$E|Y|^2<\infty$, with inner product $(Y_1, Y_2)=E({Y_1}\overline{Y_2})$; a
point in the unit sphere $S_2$ be denoted $t\equiv(1,\theta, \varphi)$; $%
X(t)=X(t,\omega)$ be mean square continuous zero mean random field on $S_2$,
$L_X^2(P(d\omega))$ be the Hilbert space generated by $X(t)$, $G$ be the
group of rotations, $G=SO(3)$.

\bigskip

1. First approach (which we present here following to \cite{Han} and \cite
{Ogu}) is based on the theory of representations of the rotation group. The
representation of the rotation group is defined by (rotational) shift
transformations on the Hilbert space generated by the random field $X(t)$.

$X(t)$, $t\in S_2$, can be regarded as a random field on $G$:

\begin{equation}  \label{Xt}
X(t)\equiv X(g_te_0), \ \ t\equiv g_t e_0, \ \ g_t\in G, \ \ g_t=g(\varphi+%
\frac{\pi}{2}, \theta, \varphi_2),
\end{equation}
where $e_0$ denotes a unit vector along the polar axis. (We suppose we have
fixed the Cartesian coordinate system with the unit vectors $%
(e_{x},e_{y},e_{z})$, with $e_{z}=e_{0}$ being the unit vector along the
polar axis of the sphere $S_{2}$ and the origin is at the center of the
sphere.)

The scalar field (\ref{Xt}) is independent of the third Euler angle or of
the rotation around $t$.

The field $X(t)$ is said to be isotropic in the wide sense, if the
correlation function $R(t_1,t_2)=E(X(t_1)\overline{X(t_2)})$ is invariant
under arbitrary rotation:
\begin{equation}  \label{RInvar}
R(t_1,t_2)=R(gt_1,gt_2), g\in G.
\end{equation}
This implies that $R(t_1,t_2)$ is a function only of $\cos\theta$, where $%
\theta$ is the angle between the unit vectors $t_1$ and $t_2$ with spherical
coordinates $(\theta_1, \varphi_1)$ and $(\theta_2, \varphi_2)$, \newline
$\cos\theta=\cos\theta_1\cos\theta_2+\sin\theta_1\sin\theta_2\cos(\varphi_1-%
\varphi_2)$.

Define the rotational transformation of the field $X(t)$:
\begin{equation}
S^{g}:X(t)\rightarrow X(g^{-1}t),\ g\in G.
\end{equation}
This transformation induces a transformation on random variables from $%
L_{X}^{2}$:
\begin{equation}
U^{g}:Y\rightarrow U^{g}Y,\ g\in G,\ Y\in L_{X}^{2}.
\end{equation}
From invariance (\ref{RInvar}) it follows that the scalar product in $%
L_{X}^{2}$, that is the covariance of $Y,Z\in L_{X}^{2}$ is invariant under $%
U^{g}$:
\begin{equation}
E(U^{g}Y\overline{U^{g}Z})=E(Y\overline{Z}),
\end{equation}
which implies that $U^{g}$ is a unitary operator on $L_{X}^{2}$. Also, we
have the group properties for $U^{g}$:
\begin{equation}
U^{g_{1}}U^{g_{2}}=U^{g_{1}g_{2}},\ \ (U^{g})^{-1}=U^{g^{-1}},\ \ U^{e}=1.
\end{equation}
This means that $U^{g}$, $g\in G$, gives the unitary representation of the
rotation group $G$ in the Hilbert space $L_{X}^{2}$. From m.s. continuity of
$X(t)$ it follows that $U^{g}$ is continuous w.r.t. $g$.

In particular, $U^g X(t)=X(g^{-1}t)$, therefore, we can write:
\begin{equation}
X(t)=X(g e_0)=U^{g^{-1}} X(e_0).
\end{equation}

Denote by $H$ the subgroup of rotations around the polar axis $e_0$ and by $%
H_t$ the subgroup of rotations around the vector $t$. Then for the scalar
filed $X(t)$ we have
\begin{equation}  \label{rotationH}
U^h X(t)=X(t), \ h\in H_t, \ U^h X(e_0)=X(e_0), \ h\in H.
\end{equation}

By the representation theory of the rotation group, the representation space
$L_{X}^{2}$ for $U^{g}$ can be decomposed into the sum of the irreducible
spaces. Correspondingly, the vector $X(e_{0})$ from $L_{X}^{2}$ can be
decomposed into the sum of the vectors of this orthogonal irreducible
spaces. Denote by $D_{l}(\Omega )$ an irreducible space of the weight $l$
representation for $U^{g}$. Denote the canonical basis for $D_{l}(\Omega )$
by
\begin{equation}
Z_{l}^{m},\ m=-l,\ldots ,l;
\end{equation}
the orthogonality relations are:
\begin{equation}
E(Z_{l}^{m}\overline{Z_{l^{\prime }}^{m^{\prime }}})=\delta _{ll^{\prime
}}\delta _{mm^{\prime }}.  \label{zlmorth}
\end{equation}

By (\ref{rotationH}), $X(e_0)$ has only 0-th canonical component for each $l$%
, that is, its irreducible decomposition can be written only in terms of $%
Z^0_l$:
\begin{equation}
X(e_0)=\sum_{l=0}^\infty \tilde{F}_l Z_l^{0},
\end{equation}
where the series converges in $L^2(P(d\omega))$, the expansion coefficients
can be given:
\begin{equation}
\tilde{F}_l =E(X(e_0) \overline{Z_l^0}),\ l=0,1,2,\ldots
\end{equation}
To obtain $X(t)$ we apply $U^{g^{-1}}$:
\begin{equation}  \label{expanmoving}
X(t)=\sum_{l=0}^\infty \tilde{F}_l U^{g^{-1}}Z_l^0.
\end{equation}

The final step is to represent $U^{g^{-1}}Z_{l}^{0}$ in terms of the
canonical basis $Z_{l}^{m}$. We know the form of the matrix of the
representation in the canonical basis, therefore, we can write:
\begin{equation}
U^{g^{-1}}Z_{l}^{0}=\sum_{m=-l}^{l}T_{0m}^{l}(g^{-1})Z_{l}^{m}=\sqrt{\frac{%
4\pi }{2l+1}}\sum_{m=-l}^{l}Y_{l}^{m}(\theta ,\varphi )Z_{l}^{m},
\end{equation}
substituting this into (\ref{expanmoving}) we obtain the spectral
representation of the field $X(t)$:
\begin{equation}
X(t)=\sum_{l=0}^{\infty }\sum_{m=-l}^{l}F_{l}Y_{l}^{m}(\theta ,\varphi
)Z_{l}^{m}  \label{reprscal}
\end{equation}
in $L^{2}(P(d\omega))$, here we have denoted $F_{l}=\tilde{F}_l\sqrt{\frac{%
4\pi}{2l+1}}$.

For the covariance function we obtain:
\begin{eqnarray}
R(\theta )&=&E(X(t_{1})\overline{X(t_{2})})=\sum_{l=0}^{\infty
}|F_{l}|^{2}\sum_{m=-l}^{l}Y_{l}^{m}(\theta _{1},\varphi _{1})\overline{%
Y_{l}^{m}(\theta _{2},\varphi _{2})}\nonumber \\&=&\frac{1}{4\pi
}\sum_{l=0}^{\infty }(2l+1)|F_{l}|^{2}P_{l}(\cos \theta ),
\label{covscal}
\end{eqnarray}
applying the addition formula for spherical harmonics
$Y_{l}^{m}(\theta ,\varphi )$
\begin{equation}
P_{l}(\cos \theta )=\frac{4\pi }{2l+1}\sum_{m=-l}^{l}Y_{l}^{m}(\theta
_{1},\varphi _{1})\overline{Y_{l}^{m}(\theta _{2},\varphi _{2})},
\end{equation}
where $P_{l}$ is the Legendre polynomial.

\begin{remark} {\rm It is interesting to note that within the
described above constructive approach for derivation of the
spectral representation we obtain the spherical harmonics
coefficients in the factorized form $F_{l}Z_{l}^{m}$ (see,
(\ref{reprscal})) with random variables $Z_{l}^{m} $ being
uncorrelated (orthonormal) in a way (\ref
{zlmorth}) and the factor $F_{l} $ is common for all $Z_{l}^{m},$ $%
m=-l,...,l.$ Note that $|F_l|^2$ is called the angular power spectrum and $%
F_lZ_l^m$ is called the random spectrum in \cite{Ogu} (or
multipole coefficients, in physical literature).}
\end{remark}

2. The spectral representation of a random field on a sphere can
be obtained as a particular case of the stochastic Peter-Weyl
theorem (see, \cite {MarinucciPeccati}, \cite{PeccatiPycke} and
the end of Section 2). Indeed, Proposition \ref{prop5} gives the
spectral decomposition for a square integrable isotropic random
field $X(g)$ on the group $SO(3)$. If we consider the restriction
of $X(g)$ on the quotient space $S_{2}=$ $SO(3)/SO(2)$, then in
the representation (\ref{spectralg}) the inner double sum will
reduce to the single sum over $m=-l,...,l,$ the functions
$D_{mn}^{l}(\varphi ,\theta ,\psi )$ will reduce to
$D_{m0}^{l}(\varphi ,\theta ,\psi )=\sqrt{\frac{4\pi
}{2l+1}}Y_{l}^{m}(\theta ,\varphi )$, and the coefficients
$a_{lmn}$ given by the formula (\ref{almn}) will simplify to
$$
a_{lmn}=\left\{
\begin{array}{ll}
0, & n=0 \\
\sqrt{2\pi }a_{lm}, & n\neq 0
\end{array}
\right.
$$
with $a_{lm}$ given by
\begin{equation}
a_{lm}=\int_{S_{2}}X(\theta ,\varphi )\overline{Y_{l}^{m}(\theta ,\varphi )}%
\sin \theta \ d\varphi \ d\theta ,  \label{alm}
\end{equation}
that is, the resulting representation appears in the form
\begin{equation}  \label{reprscal2}
X(\theta ,\varphi )=\sum_{l=0}^{\infty}\sum_{m=-l}^{l}a_{lm}Y_{l}^{m}(\theta
,\varphi ),
\end{equation}
where the random coefficients are given by (\ref{alm}).

\begin{remark} {\rm Comparing (\ref{reprscal2}) and (\ref{reprscal})
we see $a_{lm}=F_l Z_l^m$, that is we can treat the random
coefficients in
the spectral decomposition from different points of view. Note also that $%
E|a_{lm}|^2=F_l^2$ is angular power spectrum, which is more
commonly denoted as $C_l$.}\end{remark}

3. We would like to recall that the very elegant derivation of the
spectral representation of a random field on the sphere can be
given via application of classical results from analysis and
probability theory, just in three steps: (i) from the Funk-Hecke
theorem one can enjoy a complete set of eigenvalues and
orthonormal eigenfunctions for the covariance function of a mean
square continuous homogeneous isotropic random field on the
sphere, hence (ii) Mercer theorem allows to write down the
decomposition of the covariance function and finally (iii)
Karhunen theorem implies the spectral decomposition of the field
itself (for more detail see, e.g., \cite{Leonenko},
\cite{Yadrenko}).

\section{Spectral representations of vector and tensor fields}

Let now $X(t)$, $t\in S_2$, be a vector random field, that is, given a fixed
coordinate system, the filed $X(t)$ at each point $t\in S_2$ can be
represented by its coordinates $(X_1(t), X_2(t), X_3(t))$, which will change
according to the certain rules, when the coordinate system is changed. We
will also suppose that the field has zero mean.

We want to consider the field $X(t)$, which is a vector function of two
variables $(\theta, \varphi)\in S_2$, as being defined on the group of
rotations $G$, moreover we want to represent the field by a system
comprising three functions of three Euler angles $(\varphi_1, \theta,
\varphi_2)$, functions of rotations, which transform into themselves under
rotations. Thus, rotational shift transformations of each of such functions
into itself will give rise to a representation of the group $SO(3)$ and
analogously to Section 3 will lead to corresponding spectral decompositions.
To obtain such functions we will use ideas from \cite{Gel}.

Fix the Cartesian coordinate system with the unit vectors $%
(e_{x},e_{y},e_{z})$, with $e_{z}=e_{0}$ being the unit vector along the
polar axis of the sphere $S_{2}$ and the origin at the center of the sphere.
Let $(e_{r},e_{\theta },e_{\varphi })$ be the basis unit vectors for the
polar (spherical) coordinate system, this triad is associated with the point
$(1,\theta ,\varphi )$ of the unit sphere and constitutes the so-called
`moving frame'. With each rotation we can associate, however, the triad
related to an element $g$ of the rotation group. We associate with the point
$P_{0}$, the `North Pole' of the sphere, the triad $(e_{x},e_{y},e_{z})$ and
this will correspond to the identity rotation $e$. Each rotation $g$
transforms the triad at $P_{0}$ into a triad $(e_{1},e_{2},e_{3})$ and
places it at some point $P$ on the sphere, that is, the element $g$ of the
group $SO(3)$ determines the triad and the point which it is associated
with. The triad corresponding to the rotation $g=g(\varphi _{1},\theta ,0)$
is $\ (-e_{\varphi },e_{\theta },e_{r})$ and it is placed at point $%
P(\varphi ,\theta )$, where $\varphi =\frac{\pi }{2}-\varphi _{1}$, and to
the rotation $\ g=g(\varphi _{1},\theta ,\varphi _{2})$ there correspons the
triad $(e_{1},e_{2},e_{3})$ given by
\begin{eqnarray}
e_{1} &=&-e_{\varphi }\cos \varphi _{2}-e_{\theta }\sin \varphi _{2},
\label{e1} \\
e_{2} &=&-e_{\varphi }\sin \varphi _{2}+e_{\theta }\cos \varphi _{2},
\label{e2} \\
e_{3} &=&e_{r}.  \label{e3}
\end{eqnarray}
Now considering the vector $X(P)$ at point $P$ we resolve it with respect to
$(e_{1},e_{2},e_{3})$ to obtain the components $%
(X_{1}(g),X_{2}(g),X_{3}(g)). $ The component
\begin{equation}
X_{3}(g)=X_{r}(\varphi _{1},\theta ,\varphi _{2})=X_{r}(\theta ,\varphi ),\
\ \varphi =\varphi _{2}-\pi /2,  \label{normal}
\end{equation}
is normal to the surface of the sphere at point $P(\theta ,\varphi )$ and
does not depend on the third Euler angle $\varphi _{2}$, from the other two
components (lying at the tangent plane at point $P$) we form two complex
components
\begin{eqnarray*}
X_{+}(g) &=&X_{+}(\varphi _{1},\theta ,\varphi _{2})=X_{1}(g)+iX_{2}(g), \\
X_{-}(g) &=&X_{-}(\varphi _{1},\theta ,\varphi _{2})=X_{1}(g)+iX_{2}(g).
\end{eqnarray*}
These complex components can be represented as
\begin{eqnarray}
X_{+}(g) &=&e^{i\varphi _{2}}[-X_{\varphi }(\theta ,\varphi )+iX_{\theta
}(\theta ,\varphi )],  \label{xplus} \\
X_{-}(g) &=&e^{-i\varphi _{2}}[-X_{\varphi }(\theta ,\varphi )-iX_{\theta
}(\theta ,\varphi )],  \label{xminus}
\end{eqnarray}
where $X_{\varphi }(\theta ,\varphi ),$ $X_{\theta }(\theta ,\varphi )$ are
the components of the vector field in the polar coordinate system, $\varphi =%
\frac{\pi }{2}-\varphi _{1}$.

With the above three functions we have defined the vector field on the
rotation group:
$$
X(\varphi _{1},\theta ,\varphi _{2})=X(g)=(X_{+}(g),X_{-}(g),X_{r}(g)).
$$
With each rotation $g_{0}$ the normal component $X_{r}(g)=X_{r}(P)$ at point
$P(\varphi ,\theta )$ transforms into the normal component at the point $%
g_{0}^{-1}g=g_{0}^{-1}P$ \ and does not depend on the third Euler angle,
that is, on the rotation in the tangent plane.

We can consider the normal component (\ref{normal}) and the component of the
field in the tangent plane, which is represented by two functions (\ref
{xplus})-(\ref{xminus}), separately.

\begin{remark} {\rm Note that in applications vector fields on the
sphere are usually considered as vectors in the tangent plane at
each point of the sphere. Typical example of vector fields on the
sphere that arise in practice are electromagnetic fields or wind
velocity.}
\end{remark}

The notion of isotropy (in weak or second-order sense) we will define
componentwise, as isotropy property for the components (\ref{normal})-(\ref
{xminus}), that is, the covariance functions
\begin{eqnarray}
R_{+}(g_{1},g_{2}) &=&EX_{+}(g_{1})\overline{X_{+}(g_{2})},  \label{rplus} \\
R_{-}(g_{1},g_{2}) &=&EX_{-}(g_{1})\overline{X_{-}(g_{2})},  \label{rminus}
\\
R_{r}(g_{1},g_{2}) &=&EX_{r}(g_{1})\overline{X_{r}(g_{2})}  \label{rr}
\end{eqnarray}
are invariant with respect to rotations $g\in SO(3)$:
\begin{equation}  \label{covarvec}
R_{i}(g_{1},g_{2})=R_{i}(gg_{1},gg_{2}),\ \ \ \forall g,\ \ {
}i=\pm ,r,
\end{equation}
which implies $R_{i}(g_{1},g_{2})=R_{i}(g_{1}^{-1}g_{2}),$ $i=\pm ,r.$

\begin{remark} {\rm We believe that defining the notion of isotropy
in such a specific way, that is, isotropy for the components taken
in the moving basis (\ref{e1})-(\ref{e3}) is quite reasonable for
vector fields on the sphere. This will be quite analogous to the
notion of isotropic vector fields introduced in early works on
turbulence (see, e.g., \cite{Robertson}, \cite{Yaglom} and
references therein). Recall that in those works isotropy for a
vector field in strict or wide sense is defined as a property
which prescribes that probabilistic characteristics such as
probability distributions for the values of a field in a system of
points (strict sense) or just correlations (wide sense) are
invariant under translations (that is, motions, which include
rotations, shifts and reflections) of the the system of points
performed simultaneously with the same movement of the coordinate
system. In particular, for the correlation tensor of the vector
field this definition entails the following property:
$$
B(t)=G^{\ast }B(Gt)G,
$$
where $G$ is a matrix of the corresponding transformation and $G^{\ast }$ is
its transpose. (In some studies a field with the above property of its
correlation tensor is also called isotropic in a vector sense.) Considering
random fields on the sphere we are restricted to the rotations only and
define the isotropy as invariance of probabilistic characteristics under
rotations. In our approach we consider the vector field on the sphere
defined at each point of the sphere via its components in the moving basis,
which is associated with a rotation $g$. Therefore, these movements
(rotations) of the coordinate system are already incorporated into such a
representation of a vector field. And each component of the field transforms
into itself under rotation. That is why it looks quite natural to define
isotropy via the relations (\ref{covarvec}).
Moreover, this will give us the possibility to define the scalar
product, which will be invariant under rotations, in the
corresponding Hilbert spaces of random variables, and, therefore,
we can parallel all the constructions and reasonings of the
previous section, where scalar fields have been considered.}
\end{remark}

First, as we noticed above, the normal component of the field $X_{r}(\theta
,\varphi )$ does not depend on the third Euler angle and can be treated
completely in the same manner as the scalar field, and the spectral
representation of the form (\ref{reprscal}) can be written, that is, the
representation via the spherical harmonics $Y_{l}^{m}(\theta ,\varphi )$,
and the corresponding representation of the correlation function in the form
(\ref{covscal}).

Consider now the tangent component of the filed given by two functions (\ref
{xplus})-(\ref{xminus}). To obtain the spectral representation we will
follow the lines of Section 3.

Let $L_{X_{+}}^{2}$ be the Hilbert space generated by $X_{+}(g)$.

Define the rotational (or shift) transformation:
\begin{equation}
S^{g}:X_{+}(g_{0})\rightarrow X_{+}(g^{-1}g_{0}),\ g,g_{0}\in G.
\end{equation}
This transformation induces a transformation on random variables from $%
L_{X_{+}}^{2}$:
\begin{equation}
U^{g}:Y\rightarrow U^{g}Y,\ g\in G,\ Y\in L_{X_{+}}^{2}.
\end{equation}
From invariance (\ref{rplus}) it follows that the scalar product in $%
L_{X_{+}}^{2}$, that is the covariance of $Y,Z\in L_{X_{+}}^{2}$ is
invariant under $U^{g}$:
\begin{equation}
E(U^{g}Y\overline{U^{g}Z})=E(Y\overline{Z}),
\end{equation}
which implies that $U^{g}$ is a unitary operator on $L_{X_{+}}^{2}$. Also,
we have the group properties for $U^{g}$: $%
U^{g_{1}}U^{g_{2}}=U^{g_{1}g_{2}},\ \ (U^{g})^{-1}=U^{g^{-1}},\ \ U^{e}=1$.
This means that $U^{g}$, $g\in G$, gives the unitary representation of the
rotation group $G$ in the space $L_{X_{+}}^{2}$. Assuming mean square
continuity of $X_{+}(t)$, we obtain that $U^{g}$ is continuous w.r.t. $g$.

In particular, $U^{g}X_{+}(g_{0})=X_{+}(g^{-1}g_{0})$, therefore, we can
write:
\begin{equation}
X_{+}(g)=X_{+}(ge)=U^{g^{-1}}X_{+}(e),
\end{equation}
where $e$ is identity element of $G$, and
\begin{equation}
X_{+}(e)=X_{+}(e_{0})=X_{x}(e_{0})+iX_{y}(e_{0})
\end{equation}
that is, it reduces to the component of the field at the point $e_{0}=P_{0} $
(the pole of our sphere).

Denote by $H_{e_{0}}(\varphi )$ the subgroup of rotations around the polar
axis $e_{0}$ and by $H_{g}(\varphi )$ the subgroup of rotations around the
vector $ge_{0}$. Then for the filed $X_{+}(g)$ we have
\begin{equation}
U^{h}X_{+}(g)=e^{i\varphi }X_{+}(g),\ h\in H_{g}(\varphi ),\
U^{h}X_{+}(e_{0})=e^{i\varphi }X_{+}(e_{0}),\ h\in H_{e_{0}}(\varphi ).
\label{rotationHe}
\end{equation}

The representation space $L_{X_{+}}^{2}$ for $U^{g}$ can be decomposed into
the sum of the irreducible spaces. Correspondingly, the vector $X_{+}(e)$
from $L_{X_{+}}^{2}$ can be decomposed into the sum of the vectors of this
orthogonal irreducible spaces. Denote by $D_{l}(\Omega )$ an irreducible
space of the weight $l$, let
\begin{equation}
Z_{lm}^{+},\ m=-l,\ldots ,l,  \label{basisp}
\end{equation}
be the canonical basis for $D_{l}(\Omega )$, that is the basis composed from
the eigenvectors of $U^{h}$, $h\in H_{e_{0}}(\varphi )$ (rotation about $%
e_{0}$) corresponding to the eigenvalues $e^{i\varphi m}$, for these vectors
the orthogonality relations hold:
\begin{equation}
E(Z_{lm}^{+}\overline{Z_{l^{\prime }m^{\prime }}^{+}})=\delta _{ll^{\prime
}}\delta _{mm^{\prime }}.  \label{z-lm}
\end{equation}

As can be seen from (\ref{rotationHe}), $X(e_{0})$ has only `1'-st
component non-zero in the basis (\ref{basisp}), for each $l$,
therefore, its irreducible decomposition can be written only in
terms of $Z_{l1}^{+}$:
\begin{equation}
X_{+}(e_{0})=\sum_{l=0}^{\infty }F_{l}^{+}Z_{l1}^{+},
\end{equation}
where the series converges in $L^{2}(P(d\omega ))$, the expansion
coefficients can be given:
\begin{equation}
F_{l}^{+}=E(X_{+}(e_{0})\overline{Z_{l1}^{+}}),\ l=0,1,2,\ldots
\end{equation}
To obtain $X(g)$ we apply $U^{g^{-1}}$:
\begin{equation}
X_{+}(g)=\sum_{l=0}^{\infty }F_{l}^{+}U^{g^{-1}}Z_{l1}^{+}.
\label{expanmovingVec}
\end{equation}
Now it is left to notice that
\begin{equation}
U^{g^{-1}}Z_{l1}^{+}=\sum_{m=-l}^{l}T_{1m}^{l}(g^{-1})Z_{lm}^{+},
\end{equation}
substituting this into (\ref{expanmovingVec}) we obtain the spectral
representation of the field $X_{+}(g)$:
\begin{equation}
X_{+}(g)=\sum_{l=0}^{\infty
}\sum_{m=-l}^{l}F_{l}^{+}T_{1m}^{l}(g^{-1})Z_{lm}^{+}  \label{expanX+}
\end{equation}
in $L^{2}(P(d\omega ))$.

For the covariance function we obtain:
\begin{eqnarray}
R_{+}(g_{1},g_{2})&=&E(X_{+}(g_{1})\overline{X_{+}(g_{2})})=\sum_{l=0}^{\infty
}|F_{l}^{+}|^{2}\sum_{m=-l}^{l}T_{1m}^{l}(g_{1}^{-1})\overline{%
T_{1m}^{l}(g_{2}^{-1})}\nonumber \\&=&\sum_{l=0}^{\infty
}|F_{l}^{+}|^{2}T_{11}^{l}(g_{1}^{-1}g_{2}),  \label{covX+}
\end{eqnarray}
applying the addition formula (\ref{addition}) and the
orthogonality relations (\ref{z-lm}).

Analogously we obtain the spectral representation of the field $X_{-}(g)$:
\begin{equation}
X_{-}(g)=\sum_{l=0}^{\infty
}\sum_{m=-l}^{l}F_{l}^{-}T_{-1m}^{l}(g^{-1})Z_{lm}^{-}  \label{expanX-}
\end{equation}
in $L^{2}(P(d\omega ))$, where r.v. $Z_{lm}^{-}$ satisfy the orthogonality
relations $E(Z_{lm}^{-}\overline{Z_{l^{\prime }m^{\prime }}^{-}})=\delta
_{ll^{\prime }}\delta _{mm^{\prime }}$.

For the covariance function we obtain:
\begin{eqnarray}
R_{-}(g_{1},g_{2})&=&E(X_{-}(g_{1})\overline{X_{-}(g_{2})})=\sum_{l=0}^{\infty
}|F_{l}^{-}|^{2}\sum_{m=-l}^{l}T_{-1,m}^{l}(g_{1}^{-1})\overline{%
T_{-1,m}^{l}(g_{2}^{-1})}\nonumber \\ &=&\sum_{l=0}^{\infty
}|F_{l}^{-}|^{2}T_{-1,-1}^{l}(g_{1}^{-1}g_{2})  \label{covX-}
\end{eqnarray}

\bigskip Let us summarize the above reasonings in the following theorem.

\setcounter{theorem}{0}

\begin{theorem}
Let $X(t)$, $t=(\theta ,\varphi )\in S_{2}$, be a zero-mean vector random
field on the unit sphere and $(X_{\varphi }(\theta ,\varphi ),X_{\theta
}(\theta ,\varphi ),X_{r}(\theta ,\varphi ))$ be its components in spherical
coordinate system. Introducing a specific local coordinate system (\ref{e1}%
)-(\ref{e3}), which\ at every point of the sphere is related to an element
of the group of rotations, consider the representation of the field $X$ by
the components defined on the group of rotation: two complex components $%
X_{+}(g)$ and $X_{-}(g)$ given by (\ref{xplus})-(\ref{xminus}) and
the normal $X_{r}(g)=X_{r}(t)$. Suppose that these functions
$X_{+}(g),X_{-}(g)$ and $X_{r}(g)$ are weakly (second-order)
isotropic and mean square continuous. Then the functions
$X_{+}(g)$ and $X_{-}(g)$ can be expanded in the series of
generalized spherical harmonics (\ref{expanX+}) and (\ref
{expanX-}), respectively, with uncorrelated coefficients,
convergence is meant pointwise in $L^{2}(P(d\omega ))$.
Corresponding covariance functions have series representations
(\ref{covX+}) and (\ref{covX-}). The normal component can be
treated and expanded as a scalar random field.
\end{theorem}

\bigskip

Let us turn now to the tensor fields on the unit sphere. We will consider
the case of tensors of second rank.

To find the expansions for tensor fields we can apply the reasonings similar
to the above ones and represent the tensor fields by means of components
(functions on the group $SO(3)$) which transform into itself under rotations
and which will be multiplied by $e^{\pm im\varphi },$ $m=0,1,2,$ under
rotations about the axis normal to the surface of the sphere.

To obtain such components we note that components of a tensor of the second
rank are transformed under a rotation in the same way as a product of the
components of two vectors. Thus, for tensor random field on the sphere $X(t)$%
, $t\in S_{2}$, we come to the following nine components -- functions of
rotations $g=g(\varphi _{1},\theta ,\varphi _{2})$:
\begin{eqnarray}
&&X_{rr},\ \ { }X_{\varphi \varphi }+X_{\theta \theta }\pm i\left(
X_{\varphi \theta }-X_{\theta \varphi }\right) ,  \label{scal} \\
&&\left( -X_{\varphi r}-iX_{\theta r}\right) e^{i\varphi _{2}},\left(
-X_{r\varphi }-iX_{r\theta }\right) e^{i\varphi _{2}},  \label{spin1} \\
&&\left( -X_{\varphi r}+iX_{\theta r}\right) e^{-i\varphi _{2}},\left(
-X_{r\varphi }+iX_{r\theta }\right) e^{-i\varphi _{2}},  \label{spin-1} \\
&&X_{\varphi \varphi }-X_{\theta \theta }+i\left( X_{\theta \varphi
}-X_{\varphi \theta }\right) e^{2i\varphi _{2}},  \label{2i} \\
&&X_{\varphi \varphi }-X_{\theta \theta }-i\left( X_{\theta \varphi
}-X_{\varphi \theta }\right) e^{-2i\varphi _{2}}.  \label{-2i}
\end{eqnarray}

With cosmological applications in mind, we restrict our consideration to the
tangential counterpart of the field, that is, we will be interested in the
components (\ref{2i})-(\ref{-2i}) only. Moreover, we will suppose that this
tangential counterpart (the components of the tensor relative to the basis
in the tangent plane) forms a symmetric trace-free tensor, that is, $%
X_{\varphi \varphi }+X_{\theta \theta }=0$, \ $X_{\theta \varphi
}=X_{\varphi \theta }$ . In such a case the components (\ref{2i})-(\ref{-2i}%
) reduce to the form
\begin{eqnarray}
2\left( X_{\varphi \varphi }+iX_{\theta \varphi }\right) e^{2i\varphi
_{2}}&&\equiv A(g),  \label{ag} \\
2\left( X_{\varphi \varphi }-iX_{\theta \varphi }\right) e^{-2i\varphi
_{2}}&&\equiv B(g) .  \label{bg}
\end{eqnarray}

Now to derive the spectral representation for the random functions (\ref{ag}%
)-(\ref{bg}) we apply the same scheme as used above for the components (\ref
{xplus})-(\ref{xminus}) of a vector field.

We suppose that the random functions $A(g)$ and $B(g)$ are weakly (or
second-order) isotropic, that is, their covariance functions satisfy
\begin{equation}
R_{i}(g_{1},g_{2})=R_{i}(gg_{1},gg_{2}),\ \ { }\forall
g,g_{1},g_{2}\in SO(3),\ \ { }i=A,B,  \label{covar}
\end{equation}
which implies $R_{i}(g_{1},g_{2})=R_{i}(g_{1}^{-1}g_{2}),$ $i=A,B.$

Then we consider the representations of the rotation group by
rotational transformations defined on the Hilbert spaces
$L_{A}^{2}$ and $L_{B}^{2}$ generated by the random functions
$A(g)$ and $B(g),$ which will be: (1) unitary due to the
invariance of the scalar product under rotational transformations
(see, (\ref{covar})), (2)  continuous under the assumption of mean
square continuity of our random field.

We notice also that under rotations by an angle $\varphi$ about the axis
normal to the surface of the sphere (that is, around $e_{r}$) the functions $%
A(g)$ and $B(g)$ will multiply by $e^{i2\varphi }$ and $e^{-i2\varphi }$
respectively.

The representation spaces $L_{A}^{2}$ and $L_{B}^{2}$ can be decomposed into
the orthogonal sums of the irreducible spaces of weight $l$, $D_{l}^{A}$ and
$D_{l}^{B}$, introducing the canonical bases $\{Z_{lm}^{A}\}$ and $%
\{Z_{lm}^{B}\},$ $m=-l,...,l,$ in $D_{l}^{A}$ and $D_{l}^{B}$
correspondingly, we come to the representations for the functions $A(g)$ and
$B(g).$

Namely, we obtain:
\begin{equation}
A(g)=\sum_{l=0}^{\infty }\sum_{m=-l}^{l}A_{l}T_{2m}^{l}(g^{-1})Z_{lm}^{A},
\label{reprag}
\end{equation}
\begin{equation}
B(g)=\sum_{l=0}^{\infty }\sum_{m=-l}^{l}B_{l}T_{-2,m}^{l}(g^{-1})Z_{lm}^{B},
\label{reprbg}
\end{equation}
where the series converge in $L^{2}(P(d\omega ))$, r.v.'s $Z_{lm}^{A}$ and $%
Z_{lm}^{B}$\ satisfy the orthogonality relations
\begin{equation}
E(Z_{lm}^{i}\overline{Z_{l^{\prime }m^{\prime }}^{i}})=\delta
_{ll^{\prime }}\delta _{mm^{\prime }},\ \ { }i=A,B,
\end{equation}
the coefficients $A_{l}$ and $B_{l}$ are given by the formulas
\begin{equation}
A_{l}=E(A(e_{0})\overline{Z_{2l}^{A}}),\ \ { }B_{l}=E(B(e_{0})\overline{%
Z_{-2,l}^{B}}),\ l=0,1,2,\ldots
\end{equation}
For the correlation functions we obtain
\begin{equation}
R_{A}(g_{1}^{-1}g_{2})=\sum_{l=0}^{\infty
}|A_{l}|^{2}T_{22}^{l}(g_{1}^{-1}g_{2}),  \label{cova}
\end{equation}
\begin{equation}
R_{B}(g_{1}^{-1}g_{2})=\sum_{l=0}^{\infty
}|B_{l}|^{2}T_{-2,-2}^{l}(g_{1}^{-1}g_{2}).  \label{covb}
\end{equation}

We come to the following theorem.

\begin{theorem}
Let $X(t)$, $t=(\theta ,\varphi )\in S_{2}$, be a zero-mean tensor-valued
random field on the unit sphere. Suppose that the rank of the tensor is 2,
its components are given with respect to polar coordinate system and
tangential counterpart of this tensor forms a symmetric trace-free tensor.
Consider the representation of the tangential counterpart of the field $X$
by two complex components defined on the group of rotation: $\ A(g)$ and $%
B(g)$ given by (\ref{ag})-(\ref{bg}) and suppose that these
functions are weakly (second-order) isotropic and mean square
continuous. Then the functions $A(g)$ and $B(g)$ can be expanded
in the series of generalized spherical functions (\ref{reprag})
and (\ref{reprbg}), respectively, with
uncorrelated coefficients and convergence is meant pointwise in $%
L^{2}(P(d\omega ))$. Corresponding covariance functions have series
representations (\ref{cova}) and (\ref{covb}).
\end{theorem}

\setcounter{theorem}{6}
\begin{remark} {\rm We restrict ourselves here to consideration of
representations for the random functions (\ref{ag})-(\ref{bg}).
Note that the same representation holds for more general functions
(\ref{2i})-(\ref {-2i}). The functions (\ref{scal}) can be
expanded in terms of the usual
spherical harmonics $Y_{m}^{l}(\theta ,\varphi )$, the functions (\ref{spin1}%
) -- in terms of $T_{1m}^{l}$, and the functions (\ref{spin-1}) --
in terms of $T_{-1,m}^{l}$. Moreover, in the similar way the
representation of tensor fields of higher ranks then 2 can be
obtained (see also Remark \ref{rem9} below).}
\end{remark}

Let us consider more carefully the random coefficients $Z_{lm}^{+}$, $%
Z_{lm}^{-}$ and $Z_{lm}^{A}$, $Z_{lm}^{B}$ of the expansions (\ref{expanX+}%
)-(\ref{expanX-}) and (\ref{reprag})-(\ref{reprbg}). First immediate
properties of these random variables we summarize in the following lemma.

\begin{lemma}\label{lem1}
Under the condition of mean square continuity and weak isotropy of the
random function $X_{+}(g)$, the following statements hold:

(1) the random variables\ $Z_{lm}^{+}$ are uncorrelated with zero mean and
unit variance:
$$
E(Z_{lm}^{+}\overline{Z_{l^{\prime }m^{\prime }}^{+}})=\delta _{ll^{\prime
}}\delta _{mm^{\prime }};
$$

(2) for all $l$, $\overline{Z_{lm}^{+}}=(-1)^{l-m}Z_{l,-m}^{+},$ and, therefore, $%
E(Z_{lm}^{+}Z_{l^{\prime }m^{\prime }}^{+})=0$ for all $\left( l,m\right)
\neq \left( l^{\prime },m^{\prime }\right) $;

(3) for all $l$, $g\in G$, $Z_{lm}^{+}\overset{d}{=}%
\sum_{s=-l}^{l}T_{ms}^{l}(g)Z_{lm}^{+}$, in particular,
$Z_{lm}^{+}\overset{d}{=}e^{im\varphi }Z_{lm}^{+}$ $\forall
\varphi .$

The same properties hold for the random variables $Z_{lm}^{-}$ and $%
Z_{lm}^{A}$, $Z_{lm}^{B}$ under the conditions of mean square continuity and
isotropy of the corresponding fields.
\end{lemma}

\noindent{\it Proof}. All three facts follow from the very
definition of the random variables $Z_{lm}^{+}$, which have been
chosen as canonical bases in the
spaces $D_{l}$. For every $l$, $\left\{ Z_{lm}^{+},\ \ { }%
m=-l,...,l\right\} $ are orthonormal bases in orthogonal spaces $D_{l}$,
which means that (1) holds. Moreover, $\left\{ Z_{lm}^{+},\ \ { }%
m=-l,...,l\right\} $ constitute the canonical bases, for all $l$, therefore $%
U^{h}Z_{lm}^{+}=e^{im\varphi }Z_{lm}^{+}$, $h$ being the rotation around $%
e_{0}$ by an angle $\varphi $. Taking into account that $U^{g}$ is
an unitary operator for which the contravariant basis $\left\{
Z_{l}^{+(m)},\ \  m=-l,...,l\right\}$ coincides with $\left\{ \overline{Z_{lm}^{+}},\ \ { }%
m=-l,...,l\right\}$ and using the relation between the vectors of
canonical basis and its cotravariant basis:
$Z_{l}^{+(m)}=(-1)^{l-m}Z_{l,-m}^{+}$ (see \cite{Lyubarskii}), we
come to the statement (2) of the lemma. Further, we know that for
every $l$, $T^{l}(g)=\left\{ T_{ms}^{l}(g),\ \ {
}m,n=-l,...,l\right\} $ are the matrices of the representation
$U^{g}$ in the spaces $D_{l}$, with respect to the canonical basis
in $D_{l}$, that is with respect to $\left\{
Z_{lm}^{+},\ \ { }m=-l,...,l\right\} $. This implies the equality $%
U^{g}Z_{lm}^{+}=\sum_{s=-l}^{l}T_{ms}^{l}(g)Z_{lm}^{+}$. On the other hand,
isotropy of the field $X_{+}(g)$ entails the equality $U^{g}Z_{lm}^{+}%
\overset{d}{=}Z_{lm}^{+}$  $\forall l$, $\forall m$. This proves
the
statement (3) of the lemma. Since all the random variables $Z_{lm}^{+}$, $%
Z_{lm}^{-}$,$Z_{lm}^{A}$, $Z_{lm}^{B}$ have the same nature, the
properties stated in the lemma are common for all of them.

\begin{remark} {\rm Considering $2l+l$-dimensional vectors composed
from the random variables $Z_{lm}^{+}$ with fixed $l$, that is the vectors $%
Z_{l}^{+}=\left\{ Z_{lm}^{+},\ \ { }m=-l,...,l\right\} $, we can
reformulate the statement (3) of Lemma 1 in the following form: $Z_{l}^{+}%
\overset{d}{=}T^{l}(g)Z_{l}^{+}$. Statement (1) of Lemma 1 can be
reformulated in a slightly more general form as $E(Z_{l}^{+}\overline{%
Z_{l}^{+}})=I_{2l+1}$ and $E(Z_{l}^{+}\overline{Z_{l^{\prime }}^{+}})=0,$ $%
l\neq l^{\prime }$, where $I_{2l+1}$ is $\left( 2l+1\right) \times
\left( 2l+1\right) $ identity matrix and $0$ denotes $\left(
2l+1\right) \times \left( 2l^{\prime }+1\right) $ zero matrix.}
\end{remark}

The random variables $Z_{lm}^{+}$, $Z_{lm}^{-}$ and $Z_{lm}^{A}$, $%
Z_{lm}^{B} $ can be characterized further in the similar manner as
this was done in \cite{BaldiMarinucci2} for the coefficients of
expansions of scalar fields on the sphere. The properties stated
in Lemma \ref{lem1} is all that one needs to apply the same
reasonings as in \cite{BaldiMarinucci2}. Although the results
appear to be completely the same as for the scalar fields, we
present them here for completeness of the exposition and because
it is really remarkable to get the same characterization of random
coefficients within a more general construction than scalar
fields.

\begin{lemma}
Under the condition of mean square continuity and weak isotropy of the
random function $X_{+}(g)$, the the random variables $Z_{lm}^{+}$ have the
following properties: ${}{Re}Z_{lm}^{+}={}{Im}Z_{lm}^{+},$ the ratio $%
\frac{{}{Re}Z_{lm}^{+}}{{}{Im}Z_{lm}^{+}}$ is distributed
accordingly to a Cauchy distribution; ${}{Re}Z_{lm}^{+}$ and
${}{Im}Z_{lm}^{+}$
are uncorrelated with variance $\frac{1}{2}$; the marginal distributions of $%
{}{Re}Z_{lm}^{+}$ and ${}{Im}Z_{lm}^{+}$ are symmetric: ${}{Re}%
Z_{lm}^{+}=-{}{Re}Z_{lm}^{+}$, ${}{Im}Z_{lm}^{+}=-{}{Im}Z_{lm}^{+}$%
. The same properties hold for the random variables $Z_{lm}^{-}$ and $%
Z_{lm}^{A}$, $Z_{lm}^{B}$ under the conditions of mean square continuity and
isotropy of the corresponding fields.
\end{lemma}

\begin{lemma}
Let the random function $X_{+}(g)$ be mean square continuous and weakly
isotropic. If we assume in addition that $X_{+}(g)$ is Gaussian, then $%
Z_{lm}^{+},$ $l=0,1,...,$ $m=-l,...l,$ are independent Gaussian variables.
The same is true for the random variables $Z_{lm}^{-}$ and $Z_{lm}^{A}$, $%
Z_{lm}^{B}$ under the corresponding conditions on the underlying fields.
\end{lemma}

Note that to prove the independence in the above lemma the
statement (2) of Lemma \ref{lem1} is essential.

Finally, the following very important theorem can be stated.

\setcounter{theorem}{2}
\begin{theorem}
Let the random function $X_{+}(g)$ be mean square continuous and weakly
isotropic. Then for all $l$, the coefficients $\left\{ Z_{lm}^{+},\ \ { }%
m=-l,...,l\right\} $ are independent if and only if they are Gaussian. The
same is true for the random variables $Z_{lm}^{-}$ and $Z_{lm}^{A}$, $%
Z_{lm}^{B}$ under the corresponding conditions on the underlying fields.
\end{theorem}

The proof is based the Skitovich-Darmois theorem, which gives the
criterion for Gaussianity of a collection of independent random
variables via independence of some linear statistics of these
variables. For detail we refer to \cite{BaldiMarinucci2}.

Let us return to the representations (\ref{reprag}) and (\ref{reprbg}). For
the rotation $g=g(\varphi _{1},\theta ,\varphi _{2})$ we know that the
inverse rotation $g^{-1}$ is given by the Euler angles $(\pi -\varphi
_{2},\theta ,\pi -\varphi _{1}).$ Taking into account the expression for $%
T_{mn}^{l}$ (see (\ref{tlmn})) and canceling the factor $e^{\pm 2i\varphi
_{2}}$ in both sides of (\ref{reprag}) and (\ref{reprbg}) respectively, we
come to the representations
\begin{eqnarray}
2\left( X_{\varphi \varphi }(\theta ,\varphi )+iX_{\theta \varphi }(\theta
,\varphi )\right) &=&\sum_{l=0}^{\infty
}\sum_{m=-l}^{l}A_{l}T_{2m}^{l}(0,\theta ,\pi -\varphi _{1})Z_{lm}^{A} \\
&=&\sum_{l=0}^{\infty }\sum_{m=-l}^{l}A_{l}T_{2m}^{l}(0,\theta
,\pi /2+\varphi )Z_{lm}^{A},  \nonumber
\end{eqnarray}
\begin{eqnarray}
2\left( X_{\varphi \varphi }(\theta ,\varphi )-iX_{\theta \varphi }(\theta
,\varphi )\right) &=&\sum_{l=0}^{\infty
}\sum_{m=-l}^{l}B_{l}T_{-2,m}^{l}(0,\theta ,\pi -\varphi _{1})Z_{lm}^{B} \\
&=&\sum_{l=0}^{\infty }\sum_{m=-l}^{l}B_{l}T_{-2,m}^{l}(0,\theta
,\pi /2+\varphi )Z_{lm}^{B}  \nonumber
\end{eqnarray}
(we recall that $\varphi =\pi /2-\varphi _{1}$).

\bigskip

We can also reformulate the above results using Wigner
$D$-functions.

The spectral representation for the fields $A(g)$ and $B(g)$ can
be obtained as special cases of Proposition \ref{prop5}. More
precisely, due to the difference in the  convention for rotations,
we have to consider
\begin{eqnarray*}
2\left( X_{\theta \theta }+iX_{\theta \varphi }\right) e^{-2i\varphi _{2}}
&\equiv &A^{\prime }(g), \\
2\left( X_{\theta \theta }-iX_{\theta \varphi }\right) e^{2i\varphi _{2}}
&\equiv &B^{\prime }(g).
\end{eqnarray*}
Due to the special form of $A^{\prime }(g)$ and $B^{\prime }(g)$, in the
formula (\ref{spectralg}) only the functions $D_{m,-2}^{l}(g)$ and $%
D_{m,2}^{l}(g)$ will participate respectively. To use Proposition
\ref{prop5} we need also to adjust the conditions imposed on
$A^{\prime }(g)$ and $B^{\prime }(g).$ Namely, we demand the
fields $A^{\prime }(g)$ and $B^{\prime }(g)$ to be square
integrable and strictly isotropic (as described in Proposition
\ref{prop5}). Then we can write the decompositions
\begin{equation}
A^{\prime }(g)=\sum_{l=0}^{\infty }\sum_{m=-l}^{l}a_{lm}\sqrt{\frac{2l+1}{%
8\pi ^{2}}}D_{m,-2}^{l}(g),  \label{a'g}
\end{equation}
\begin{equation}
B^{\prime }(g)=\sum_{l=0}^{\infty }\sum_{m=-l}^{l}b_{lm}\sqrt{\frac{2l+1}{%
8\pi ^{2}}}D_{m,2}^{l}(g),  \label{b'g}
\end{equation}
where
\begin{equation}
a_{lm}=\int_{G}A^{\prime }(g)\sqrt{\frac{2l+1}{8\pi ^{2}}}\overline{%
D_{m,-2}^{l}}(g)dg,
\end{equation}
\begin{equation}
b_{lm}=\int_{G}B^{\prime }(g)\sqrt{\frac{2l+1}{8\pi ^{2}}}\overline{%
D_{m,2}^{l}}(g)dg,
\end{equation}
the convergence of series is in $L_{2}(P(d\omega )\times dg)$ and poinwise
in $L_{2}(P(d\omega ))$.

\bigskip

We can also deduce from (\ref{a'g})-(\ref{b'g}) the representations
\begin{equation}
2\left( X_{\theta \theta }(\theta ,\varphi )+iX_{\theta \varphi }(\theta
,\varphi )\right) =\sum_{l=0}^{\infty }\sum_{m=-l}^{l}a_{lm}\sqrt{\frac{2l+1%
}{8\pi ^{2}}}D_{m,-2}^{l}(\varphi ,\theta ,0),  \label{spin-2}
\end{equation}
\begin{equation}
2\left( X_{\theta \theta }(\theta ,\varphi )-iX_{\theta \varphi }(\theta
,\varphi )\right) =\sum_{l=0}^{\infty }\sum_{m=-l}^{l}b_{lm}\sqrt{\frac{2l+1%
}{8\pi ^{2}}}D_{m,2}^{l}(\varphi ,\theta ,0).  \label{spin2}
\end{equation}
\setcounter{theorem}{8}
\begin{remark}\label{rem9} {\rm We can see that the key step to
derive the spectral decomposition for vector and tensor fields on
the sphere is to consider the components which transform in a
specific way under rotations about the axis $e_{r}$. In harmonic
analysis on the sphere there exists a special class of functions
which are defined relatively their behavior under the rotations of
the basis vectors in the plane tangent to the sphere (i.e.
rotations around $e_{r}$). These are so-called spin-weighted
functions. We recall the definition of these functions (see,
\cite{Goldberg}, \cite{Newman}, \cite{Thorne}), however here we
use different convention for rotation than the above authors to be
consistent with the literature on CMB polarization. First of all,
we note that these functions, similarly to the vector- or
tensor-valued functions, are defined relatively to some coordinate
system. Let at any point of the sphere we have defined three
orthogonal vectors: one
normal (radial) $e_{r}$ and two tangential, $(e_{1},e_{2})$. A function $%
f(\theta ,\varphi )$ defined on the sphere $S_{2}$ is said to have spin $s$
if under a right-handed rotations of basis vectors $(e_{1},e_{2})$ by an
angle $\psi $ it transforms as $f^{\prime }(\theta ,\varphi )=e^{-is\psi
}f(\theta ,\varphi )$. Spin functions can be equivalently defined as
evaluation at $\varphi _{2}=0$ of any function in $L_{2}(SO(3))$ resulting
from an expansion for fixed index $n$ in the Wigner $D$-functions $%
D_{mn}^{l}(\varphi _{1},\theta ,\varphi _{2})$. Thus, the functions $%
D_{mn}^{l}(\varphi ,\theta ,0)$ or $D_{m,-n}^{l\ast }(\varphi ,\theta ,0)$
define an orthogonal basis for the expansions of spin $n$ functions in $%
L_{2}(S_{2})$. After normalization in $L_{2}(S_{2})$, these basis functions,
which are called the spin-weighted spherical harmonics of spin $n$, are
given in a factorized form in terms of the real Wigner $d$-functions $%
d_{mn}^{l}(\theta )$ and the complex exponentials $e^{im\varphi }$ as
follows:
\begin{equation}
_{n}Y_{m}^{l}(\theta ,\varphi )=\sqrt{\frac{2l+1}{4\pi }}d_{m,-n}^{l}(\cos
\theta )e^{im\varphi }.
\end{equation}

We refer for the rigorous mathematical theory for spin-weighted functions to
\cite{GellerMarinucci2008}, \cite{GellerLanMarinucci2009}.

For example, $X_{\theta \theta }(\theta ,\varphi )\pm iX_{\theta
\varphi }(\theta ,\varphi)$ are spin $\pm 2$ functions and the
decompositions (\ref {spin-2})-(\ref{spin2}) are decompositions in
spin $\pm 2$ spherical harmonics.}
\end{remark}
\noindent \textbf{Example.} In cosmological studies the observable
linear polarization field is described in terms of two Stokes'
parameters $Q$ and $U $ defined in a given direction with respect
to the particular choice of axes in a plane perpendicular to the
direction of observation. When reference frame is rotated around
the direction of observation (direction of propagation) $Q$ and
$U$ transform like the components of a 2-dimensional 2-d rank
symmetric trace-free tensor. Therefore, linear polarization can be
conveniently described as a tensor-valued field on a sphere
$$
\frac{1}{2}\left(
\begin{array}{cc}
Q & U \\
U & -Q
\end{array}
\right)
$$
defined with respect to the polar basis vectors $(e_{\theta
},e_{\varphi })$, or, equivalently, via combinations $Q\pm iU$.

Using the above theory, we can expand $Q(\theta, \varphi)\pm
iU(\theta, \varphi)$ in the spin $\pm 2$ spherical harmonics. For
this we must impose the assumption of strict isotropy, as in
Proposition \ref{prop5}, then we can write
\begin{equation}\label{pol}
P(\theta, \varphi)\equiv Q(\theta, \varphi)+iU(\theta ,\varphi )
=\sum_{l=0}^{\infty }\sum_{m=-l}^{l} {_2}p_{lm}\,\, {_2}Y_{m}^{l}(\theta
,\varphi),
\end{equation}
\begin{equation}\label{pol*}
P^*(\theta, \varphi)\equiv Q(\theta, \varphi)-iU(\theta ,\varphi )
=\sum_{l=0}^{\infty }\sum_{m=-l}^{l} {_{-2}}p_{lm}\,\, {_{-2}}%
Y_{m}^{l}(\theta ,\varphi)
\end{equation}
both in $L_{2}(P(d\omega )\times dg)$ and pointwise in
$L_{2}(P(d\omega ))$.

However, as our first approach suggests, the above decompositions
are also hold for mean square continuous random fields under the
assumption on weak (second-order) isotropy pointwise in
$L_{2}(P(d\omega ))$. From the  theory presented in this section
we have some characteristic properties of random coefficients in
the above series (for mean square continuous weakly isotropic
random fields). We know, in particular, that the coefficients are
uncorrelated, for Gaussian fields the coefficients are Gaussian
independent variables, and the coefficients are independent if and
only if they are Gaussian. Thus, with the rigorous probabilistic
framework provided in the paper, expansions of the form
(\ref{pol})-(\ref{pol*}), which are widely used in the literature
on CMB polarization, now become operational.

\end{document}